\documentclass[12pt]{amsart}
\let\oldlabel=\label
\def\prellabel{\marginparsep=1em\marginparwidth=44pt
    \def\label##1{\oldlabel{##1}\ifmmode\else\ifinner\else
         \marginpar{{\footnotesize\ \\\normalshape\tt 
                    ##1}}\fi\fi}}
\def\Q{\ifhmode\textqed\fi
   \ifmmode\ifinner\quad\Qsymbol\else\dispqed\fi\fi}
\def\textqed{\unskip\nobreak\penalty50
    \hskip2em\hbox{}\nobreak\hfil\Qsymbol
    \parfillskip=0pt \finalhyphendemerits=0}
\def\dispqed{\rlap{\qquad\Qsymbol}}

\usepackage{amscd}

\def\and{\operatorname{and}}
\def\fa{\operatorname{for~ all}}
\def\ker{\operatorname{ker}}
\def\im{\operatorname{im}}
\def\dim{\operatorname{dim}}
\def\depth{\operatorname{depth}}

\def\mm{{\mathfrak m}}

\def\CC{{\mathbb C}}

\def\ZZ{{\mathbb Z}}

\def\mm{{\frak m}}
\def\lrar{{\longrightarrow}}
\def\R{{\mathcal R}}
\def\G{{\mathcal G}}
\def\wrt{{ with respect to }}

\def\rr{{ Ratliff-Rush }}

\def\lm{{\lambda}}
\def\bl{\begin{lemma}}
\def\el{\end{lemma}}
\def\bt{\begin{theorem}}
\def\et{\end{theorem}}
\def\ben{\begin{enumerate}}
\def\een{\end{enumerate}}
\def\bpf{\begin{proof}}
\def\epf{\end{proof}}
\def\beqn{\begin{eqnarray*}}
\def\eeqn{\end{eqnarray*}}
\def\bd{\begin{definition}}
\def\ed{\end{definition}}
\def\bp{\begin{proposition}}
\def\ep{\end{proposition}}
\def\bc{\begin{corollary}}
\def\ec{\end{corollary}}

\def\rr{{ Ratliff-Rush }}

\let\epsilon=\varepsilon
\let\phi=\varphi
\let\theta=\vartheta

\newtheorem{lemma}{Lemma}[section]
\newtheorem{corollary}[lemma]{Corollary}
\newtheorem{theorem}[lemma]{Theorem}
\newtheorem{proposition}[lemma]{Proposition}

\theoremstyle{definition}
\newtheorem{definition}[lemma]{Definition}

\newtheorem{example}[lemma]{Example}

\advance\headheight1.15pt

\begin{document}
\title[Grothendieck-Serre formula]{Grothendieck-Serre formula and bigraded
Cohen-Macaulay Rees algebras}

\author{A. V. Jayanthan \and J. K. Verma}
\thanks{The
first author is supported by the National Board for Higher Mathematics, India}
\address{Department of Mathematics, Indian Institute of Technology
Bombay, Powai, Mumbai, India - 400076}
\email{jayan@math.iitb.ac.in}
\email{jkv@math.iitb.ac.in}

\maketitle
{\it Dedicated to Prof. Dr. J\"urgen Herzog on the occasion of his sixtieth
birthday}
\begin{abstract} The Grothendieck-Serre formula for the difference between
the Hilbert function and Hilbert polynomial of a graded algebra is
generalized for bigraded standard algebras. This is used to get a similar
formula for the difference between the Bhattacharya function and
Bhattacharya polynomial of two $\mm$-primary ideals $I$ and $J$ in a local
ring $(A,\mm)$ in terms of local cohomology modules of Rees algebras of
$I$ and $J.$ The cohomology of a  variation of the Kirby-Mehran complex
for bigraded Rees algebras   is studied which is used to characterize the
Cohen-Macaulay property of bigraded Rees algebra of  $I$ and $J$ for two
dimensional Cohen-Macaulay local rings. 

\end{abstract}
\setcounter{page}{1}
\tableofcontents
\pagebreak
\section{Introduction}

Let $R = \oplus_{n \geq 0}R_n$ be a finitely generated standard graded
algebra over an Artinian local ring $R_0$. Let $\lambda$ denote length. 
The Hilbert function of $R$,
$H(R,n)=\lambda_{R_0}(R_n)$, is given by a 
polynomial $P(R,n) $ for $n \gg 0.$ The Grothendieck-Serre
formula expresses the difference $H(R,n)-P(R,n)$   in terms of lengths of
graded components of the local cohomology modules of $R$ with support in
the irrelevant ideal  
$R_{+} = \oplus_{n > 0 }R_n$ of $R.$ 
We shall prove a version of this formula  in Section 2 for bigraded
standard algebras over Artinian local rings. We need this generalization to
find necessary and sufficient conditions for the Cohen-Macaulay property
of bigraded Rees algebras. These conditions involve the coefficients of
the Bhattacharya polynomial of two $\mm$-primary ideals in a local ring
$(R,\mm).$ 

    To be more precise, let $I$ and $J$ be
$\mm$-primary ideals in a $d$-dimensional local ring $(R,\mm).$ The function 
$B(r,s)=\lambda(R/I^rJ^s)$ is called the Bhattacharya function of $I$ and
$J$ \cite{b}. Bhattacharya proved in \cite{b} that this function
is given by a polynomial $P(r,s)$ for $r, s \gg 0$. We represent the Bhattacharya
polynomial $P(r,s)$ corresponding to $B(r,s)$ by 
$$
P(r, s) = \sum_{i+j \leq d} e_{ij}{r \choose i}{s \choose j}
$$ 
where  $e_{ij} \in \ZZ.$  The integers $e_{ij}$ for which $i+j = d$ were
termed as mixed multiplicities of $I$ and $J$ by Teissier and Risler
in \cite{t}. We write $e_j(I|J)$ for $e_{ij}$ when $i+j = d$. 

The bigraded version of the Grothendieck-Serre formula, proved in
Section 2, allows us to
express the difference of the Bhattacharya function and  Bhattacharya  
polynomial of two $\mm$-primary ideals $I$ and $J$  in terms of lengths of
bigraded components of local cohomology
modules of the extended Rees algebra of  $I$ and $J.$ This is done in section
5 of the paper.  

 In Section 3 we prove some preliminary 
results about Ratliff-Rush closure of products of ideals. In Section 4 we
present a variation on a complex first defined by Kirby 
and Mehran in \cite{km}. The cohomology of this complex is related to
the local cohomology of Rees algebras of two ideals. An analysis of this 
relationship yields a formula for the constant term  of the 
Bhattacharya polynomial $P(r,s)$. This formula is used to prove the
characterization of Cohen-Macaulay property of bigraded Rees algebras
mentioned above.

\vspace{.5cm}

\noindent
{\bf Acknowledgment :} We thank the referee for a careful reading,
suggesting several improvements and pointing out related references. 

\section{Grothendieck-Serre difference formula for bigraded algebras}
We begin by establishing the notation for bigraded algebras. A ring
$A$ is called a bigraded algebra if $A = \oplus_{r, s \in
\ZZ}A_{(r,s)}$ where each $A_{(r,s)}$ is an additive subgroup of $A$
such that $A_{(r,s)}\cdot A_{(l,m)} \subseteq A_{(r+l, s+m)} \fa
\;\; (r,s), (l,m) \in \ZZ^2$. We say that $A$ is a standard bigraded
algebra if $A$ is finitely generated, as an $A_{(0,0)}$-algebra, by 
elements of degree $(1,0)$ and $(0,1).$ The elements of $A_{(r,s)}$ 
are called bihomogeneous of degree $(r,s)$.
An ideal $I$ of $A$ is said to be bihomogeneous if $I$ is generated by
bihomogeneous elements.  The ideal of $A$ generated by
elements of degree $(r,s)$, where $r + s \geq 1$ is denoted by $A_{+}$
and the ideal generated by elements of degree $(r,s)$, where 
$r,s \geq 1$ is denoted by $A_{++}$. An $A$-module $M$ is called bigraded
if $M
= \oplus_{r,s \in \ZZ}M_{(r,s)}$, where $M_{(r,s)}$ are additive
subgroups of $M$ satisfying $A_{(r,s)}\cdot M_{(l,m)} \subseteq
M_{(r+l, s+m)}$ for all $r,s,l,m \in \ZZ$. It is known that 
when $A_{(0,0)}$ is Artinian and $M$
is a finitely generated bigraded $A$-module, the function
$\lambda_{A_{(0,0)}}(M_{(r,s)})$, called Hilbert function of $M$, is 
finite for all $r,s$ and coincides with a polynomial for $r, s \gg 0$.
In this section we express the difference between the Hilbert function
and the Hilbert polynomial in terms of the Euler characteristic of local
cohomology modules. For an ideal $I$ in $A$ and an $A$-module $M,$ let
$H^i_I(M)$ denote the $i$-th local cohomology module of $M$ with
respect to $I$. We refer the reader to [BS] for properties of
local cohomology modules. Note that when $I$ is a bihomogeneous ideal
in a bigraded algebra $A$ and $M$ is a bigraded $A$-module, the local 
cohomology modules $H^i_I(M)$ have
a natural bigraded structure inherited from $A$ and $M$.

Throughout this section $(A,\mm)$ will denote a $d$-dimensional 
Noetherian local ring unless stated otherwise. Let $X = (X_1,
\ldots, X_m)$ and $Y = (Y_1, \ldots, Y_n)$ be two sets of
indeterminates. Let $R= A[X_1, \ldots, X_m, Y_1, \ldots, Y_n]$. 
We assign the grading deg $X_i = (1,0)$ for $i = 1, \ldots, m$ and 
deg $Y_i = (0,1)$ for $i = 1, \ldots, n$ so that
$R$ is a standard bigraded algebra. We write $R_{(r,s)}$ for the
$A$-module generated by products of monomials of degree $r$ in $X$ and
degree $s$ in $Y$. In the next lemma we establish finite generation
over $A$ of the bigraded components of the local cohomology modules of
$R$ \wrt $X$ and $Y$ respectively. 

The results in this section are not
new. They are folklore in the multigraded case. Lemma \ref{pollem}
follows from Lemma 2.2 and Corollary 2.3 of \cite{cht} when $A$ is a
field. Theorem \ref{imp2} and Theorem \ref{gro} follow from Lemma 4.2
and Lemma 4.3 of \cite{kt}. We refer the reader to Lemma 2.1 of 
\cite{o2}, Theorem 9.1 of \cite{sn} and Section 1 of \cite{k}.

Although the results in the section are not new, we have provided easy
proofs so that these results are accessible to readers not familiar
with sheaf cohomology.

\bl\label{pollem}
Let $R = A[X_1, \ldots, X_m, Y_1, \ldots, Y_n].$ Then

\begin{enumerate}

\item[(i)] $H^i_X(R) = 0$ for all $i \neq m$ and 
$H^i_Y(R) = 0$ for all $i \neq n.$

\item[(ii)]  $H^m_{X}(R)_{(r,s)} = 0$ for all $r > -m$ and,
$H^n_{Y}(R)_{(r,s)} = 0$ for all $s > -n.$

\item[(iii)] $H^m_{X}(R)_{(r,s)} \mbox{ and } H^n_{Y}(R)_{(r,s)} $ are
finitely generated A-modules for all $r,s \in \ZZ.$

\end{enumerate}
\el

\bpf 
(i) is standard. 

(ii). Induct on $m$.
Let $m =0$. Then $H^0_{(0)}(R) = R = A[Y].$ Therefore
$H^0_{(0)}(R)_{(r,s)} = 0$ for all $r > 0.$ Suppose $m > 0.$ 
Let $\bar{R} = R/X_mR$ and $(\bar{X}) = (X_1, \ldots, X_{m-1})$.
Consider the short exact sequence
$$
0 \lrar R(-1,0) \buildrel{.X_m}\over\longrightarrow R \lrar
\bar{R} \lrar 0.
$$
By the change of ring principle, $H^i_{(X)}(\bar{R}) =
H^i_{(\bar{X})}(\bar{R})$. Since $(\bar{X})$ is generated by $m-1$
indeterminates, $H^i_{(\bar{X})}(\bar{R}) = 0$ for all $i \neq m-1$. 
Therefore we get the following long exact sequence 
\begin{eqnarray}\label{topseq}
0 & \lrar & H^{m-1}_{(\bar{X})}(\bar{R}) \lrar
H^m_{(X)}(R)(-1,0) \buildrel{.X_m}\over\longrightarrow H^m_{(X)}(R) \lrar 0.
\end{eqnarray}
By induction hypothesis, for all $r > -m+1,~
H^{m-1}_{(\bar{X})}(\bar{R})_{(r,s)} = 0$. Hence for $r > -m+1$ 
we get an exact sequence 
$$
0 \lrar H^m_{(X)}(R)_{(r-1,s)} \buildrel{.X_m}\over\longrightarrow 
H^m_{(X)}(R)_{(r,s)} \lrar 0.
$$
Let $z \in H^m_{(X)}(R)_{(r-1,s)}$. Pick the smallest $l \geq 1$,
such that $X_m^lz = 0$. Then $X_m(zX_m^{l-1}) = 0$. Therefore $z = 0$.
Hence $H^m_{(X)}(R)_{(r,s)} = 0$ for all $r > -m$. Similarly one can
can show that $H^n_{Y}(R)_{(r,s)} = 0$ for all $s > -n.$ 

(iii) We need to show that
$H^m_{(X)}(R)_{(r,s)}$ is finitely generated  for all $r \leq -m$. Apply induction on
$m$. It is clear for $m = 0$. Assume the statement for $m-1$.
Now apply decreasing induction on $r$.
When $ r = -m+1$, $H^{m-1}_{(\bar{X})}(\bar{R})_{(-m+1,s)} \cong 
H^m_{(X)}(R)_{(-m,s)}$, by (\ref{topseq}) and (ii).
By induction hypothesis on $m$,
$H^{m-1}_{(\bar{X})}(\bar{R})_{(-m+1,s)}$ is finitely generated
hence so is $H^m_{(X)}(R)_{(-m,s)}$. Now for $r < -m+1$
we have the short exact sequence
$$
0 \lrar H^{m-1}_{(\bar{X})}(\bar{R})_{(r,s)} \lrar
H^m_{(X)}(R)_{(r-1,s)} \buildrel{.X_m}\over\longrightarrow
H^m_{(X)}(R)_{(r,s)} \lrar 0.
$$
By induction on $r$, $H^m_{(X)}(R)_{(r,s)}$ is finitely generated
and $H^{m-1}_{(\bar{X})}(\bar{R})_{(r,s)}$ is finitely generated 
by induction on $m$. Therefore  
$H^m_{(X)}(R)_{(r-1,s)}$ is finitely generated. Similarly $H_Y^n(R)_{(r,s)}$ is
finitely generated for all $r,s \in \ZZ$.
\epf

\bl\label{imp1}
\ben
\item[(i)]   $H^i_{R_{++}}(R) = 0$ for all $i \neq m, n$ and $m+n-1.$
\item[(ii)]  $H^i_{R_{++}}(R)_{(r,s)} = 0$ for $r,s \gg 0$ and $i \geq 0.$
\item[(iii)] $H^i_{R_{++}}(R)_{(r,s)}$ is a finitely generated $A$-module for all $i
\geq 0$ and $r,s \in \ZZ.$
\een
\el

\bpf
First note that $R_{++} = (X_1, \ldots, X_m) \cap (Y_1, \ldots, Y_n)$.
Set $X = (X_1, \ldots, X_m)$,  $Y = (Y_1, \ldots, Y_n)$ and $R_{+} =
X + Y$. Consider the Mayer-Vietoris sequence :
\begin{eqnarray}\label{1}
\cdots \lrar H^i_{R_{+}}(R) \lrar H^i_X(R) \oplus H^i_Y(R) \lrar
H^i_{R_{++}}(R) \lrar H^{i+1}_{R_{+}}(R) \lrar \cdots 
\end{eqnarray}

(i). If $i \neq m, n, m+n-1, ~ H^i_X(R) = H^i_Y(R) = H^{i+1}_{R_{+}}(R)
=0$. Hence $H^i_{R_{++}}(R) = 0$ for $i \neq m, n, m+n-1$. 

\noindent
By Theorem 2.2.4 of \cite{bl} and Lemma \ref{pollem}, (ii) and (iii)
are satisfied by $H^i_X(R), H^i_Y(R)$ and $H^{i+1}_{R_{+}}(R)$. Hence
(ii) and (iii) are satisfied by $H^i_{R_{++}}(R)$. 
\epf

\bt\label{imp2}
Let $R = \oplus_{r,s \geq 0}R_{(r,s)}$ be a finitely generated 
standard bigraded 
algebra over a Noetherian local ring $R_{00} = (A,\mm)$. Let 
$M$ be a finitely generated  bigraded $R$-module. Then

\ben

\item[(i)]  $H^i_{R_{++}}(M)_{(r,s)} = 0$ for all $r, s \gg 0$ and $i
\geq 0.$

\item[(ii)] $H^i_{R_{++}}(M)_{(r,s)}$ is a finitely generated $A$-module for all
$r,s \in \ZZ$ and $i \geq 0$.

\een
\et

\bpf
As $R$ is standard bigraded $R \cong A[X_1, \ldots, X_m, Y_1, \ldots, Y_n]/I$
for a bihomogeneous ideal $I$.
Consider $M$ as a bigraded $S = A[X_1, \ldots, X_m, Y_1, \ldots,  Y_n]$-module. Then by the change of ring principle
$H^i_{R_{++}}(M) = H^i_{S_{++}}(M)$ for all $i \geq 0$. Therefore,
without loss of generality, 
we may assume that $R = A[X_1, \ldots, X_m, Y_1, \ldots, Y_n]$.
Since $M$ is a finitely generated bigraded $R$-module, there exists
a free $R$-module $F = \oplus_{j=1}^sR(m_j), m_j \in \ZZ^2$ 
and a short exact sequence of finitely generated bigraded $R$-modules
$$
0 \lrar K \lrar F \lrar M \lrar 0.
$$
Consider the corresponding long exact sequence of local cohomology
modules
$$
\cdots \lrar H^i_{R_{++}}(K) \lrar H^i_{R_{++}}(F) \lrar
H^i_{R_{++}}(M) \lrar H^{i+1}_{R_{++}}(K) \lrar \cdots 
$$
By  Lemma \ref{imp1},  (i) and (ii) are true for
$H^i_{R_{++}}(F)$. We prove the theorem by decreasing induction on
$i$. Since $H^i_{R_{++}}(M) = 0$ for $i \gg 0$, (i) and (ii) obviously
hold for $i \gg 0$. By induction $H^{i+1}_{R_{++}}(K)$ has properties
(i) and (ii). Hence $H^i_{R_{++}}(M)$ satisfies (i) and
(ii).
\epf

\bt\label{gro}
Let $R = \oplus_{r,s \geq 0}R_{(r,s)}$ be a finitely generated standard
bigraded algebra with $R_{00} = (A, \mm)$, an Artinian local ring and 
let $M = \oplus_{r ,s \geq 0}M_{(r,s)}$ be a bigraded 
finite $R$-module. Put $B_M(r,s) = \lambda_A(M_{(r,s)}).$ Let
$P_M(r,s)$ denote the Hilbert polynomial corresponding to the function
$B_M(r,s)$.  Then for all $r,s \in \ZZ$,
$$
B_M(r,s) - P_M(r,s) = \sum_{i\geq0}(-1)^i\lambda_{A}(H^i_{R_{++}}(M)_{(r,s)}).
$$
\et

\begin{proof}
Write $R = A[x_1, \ldots, x_m, y_1, \ldots, y_n]$ with deg $x_i$ =
(1,0) and deg $y_i$ = (0,1). We prove the theorem by induction on
$m+n$. Suppose $m+n = 0$. Then $M_{(r,s)} = 0$ for $r,s \gg 0$. Hence
$P_M(r,s) = 0$. Since $\dim M$ = 0, we have $H^i_{R_{++}}(M) = 0$ for
all $i > 0$ and $H^0_{R_{++}}(M) = M$. Therefore $B_M(r,s) =
\lambda_A(H^0_{R_{++}}(M)_{(r,s)})$.
\\ \\ 
Now suppose $m+n > 0$. If $m = 0$ or $n = 0$, the result reduces to
Theorem 2.2.2 of \cite{bl}. 
Let $m > 0$ and $n > 0$. Consider the exact sequence of finitely
generated
bigraded $R$-modules
\begin{eqnarray}\label{sesgro}
0 \lrar K \lrar M(-1,0)
\buildrel{.x_m}\over\longrightarrow M \lrar C \lrar 0.
\end{eqnarray}
For any finitely generated bigraded $R$-module $N$, define
$$
\chi_N(r,s) = \sum_{i \geq 0}(-1)^i\lambda_A(H^i_{R_{++}}(N)_{(r,s)})
$$
and 
$$
f_N(r,s) = B_N(r,s) - P_N(r,s).
$$
Since $H^i_{R_{++}}(N(-\mu,0))_{(r,s)} = H^i_{R_{++}}(N)_{(r-\mu,s)}$, 
it follows that $\chi_{N(-\mu,0)}(r,s) = \chi_N(r-\mu,s)$. 
Thus from (\ref{sesgro}), we get 

$$\chi_M(r-1,s) - \chi_M(r,s) = \chi_K(r,s) - \chi_C(r,s)$$ 
and 
$$f_M(r-1,s) - f_M(r,s) = f_K(r,s) - f_C(r,s)$$
for all $r,s \in \ZZ$. Let $\bar{R} = R/x_mR \cong A[\bar{x}_1,
\ldots, \bar{x}_{m-1}, \bar{y}_1, \ldots, \bar{y}_n]$. Since $x_mK = 0
= x_mC$, we can consider $K$ and $C$ as $\bar{R}$-modules. By the change
of ring principle, 
$$H^i_{R_{++}}(K) \cong H^i_{\bar{R}_{++}}(K)
\mbox{ and }  H^i_{R_{++}}(C) \cong H^i_{\bar{R}_{++}}(C)$$ for all $i \geq
0$. By induction $f_K(r,s) = \chi_K(r,s)$ and $f_C(r,s) =
\chi_C(r,s)$. Therefore we have $\chi_M(r,s) - \chi_M(r-1, s) =
f_M(r,s) - f_M(r-1, s)$ for all $(r,s) \in \ZZ^2$. Consider
the exact sequence (\ref{sesgro}) with the map, multiplication by $y_n$.
Proceeding as in the above case we get that
$\chi_M(r,s) - \chi_M(r, s-1) = f_M(r,s) - f_M(r,s-1)$.
By Theorem \ref{imp2}, $\chi_M(r,s) = 0$ for $r,s \gg 0$ and
clearly $f_M(r,s) = 0$ for $r,s \gg 0$.
Set $h = \chi_M - f_M$. then
$h(r,s) = 0$ for all $r,s \gg 0$ and we have $h(r,s) = h(r-1,s),
h(r,s) = h(r, s-1)$ for all $r, s$. Therefore $h = 0$ and
$$
B_M(r,s) - P_M(r,s) =
\sum_{i\geq0}(-1)^i\lambda_{A}(H^i_{R_{++}}(M)_{(r,s)}).
$$
\end{proof}

\section{Ratliff-Rush closure of products of ideals}

Let $A$ be a commutative ring and $K \subset I$ be ideals of $A$. We 
say that $K$ is a {\it reduction} of $I$ if there exists an integer $r 
\geq 1$ such that $I^r = KI^{r-1}$. 
The smallest integer $r$ satisfying this equation is
called the reduction number, $r_K(I)$, of $I$ \wrt $K.$
We say that $K$ is a minimal reduction of $I$ if $K$ is minimal \wrt
inclusion among all reductions of $I$. We refer the reader to
\cite{nr}
for basic facts about reductions of ideals.

Let $(A, \mm)$ be a local ring and $I$ be an ideal of $A$. The stable
value of the sequence $\{I^{n+1} : I^n\}$ is called the Ratliff-Rush
closure of $I$, denoted by $\tilde{I}$. An ideal $I$ is said to be
Ratliff-Rush if $\tilde{I} = I$. In this section we discuss the
concept of the Ratliff-Rush closure for the product of two ideals. 

The following proposition summarizes some basic properties of \rr
closure found in [RR].
\begin{proposition} Let $I$ be an ideal containing regular element in
a Noetherian ring $A$. Then
\begin{enumerate}
\item $I \subseteq \tilde{I}$ and $\widetilde{(\tilde{I})} = \tilde{I}$.  
\item $(\tilde{I})^n = I^n$ for $n \gg 0$. Hence if $I$ is 
$\mm$-primary, the Hilbert polynomial of $I$ and $\tilde{I}$ are same.
\item $\tilde{(I^n)} = I^n$ for $n \gg 0$.
\item If $(x_1, \ldots, x_g)$ is a minimal reduction of $I$, then 
$\tilde{I} = \cup_{n \geq 0} I^{n+1} : (x_1^n, \ldots, x_g^n)$.
\end{enumerate}
\end{proposition}

We show that the Ratliff-Rush closure for product of two ideals can be
computed from complete reductions, a generalization of reductions of
ideals introduced by Rees in \cite{r2}. 

Let $(A,\mm)$ be a $d$-dimensional local ring. Let $I_1, \ldots, I_r$ 
be $\mm$-primary ideals of $(A,\mm)$. Let $\left(x_{ij}\right)$ with 
$x_{ij} \in I_i,$ for all $j = 1, \ldots, d \mbox{ and } i = 1, 
\ldots, r$, be a system of elements in $A$. Put $y_j = x_{1j}x_{2j}
\ldots x_{rj}, j = 1, \ldots, d$. Then the system 
of elements $(x_{ij})$ is said to be a {\it complete reduction} 
of the sequence of ideals $I_1, \ldots, I_r$ if $(y_1, \ldots, y_d)$ 
is a 
reduction of $I_1\ldots I_r$. In \cite{r2} Rees proved the existence of 
complete reductions when the residue field of $A$ is infinite. 

\bl\label{rrprod}
Let $I$ and $J$ be ideals of $A$.  Then we have 
\begin{enumerate}
\item[(i)] $\widetilde{IJ} = \displaystyle{\bigcup_{r,s \geq 0}}
I^{r+1}J^{s+1} : I^rJ^s$.

\item[(ii)] $(\widetilde{I^aJ^b}) = \displaystyle{\bigcup_{k \geq 0}}I^{a+k}J^{b+k}
: I^kJ^k$.
\item[(iii)] If $I$ and $J$ are $\mm$-primary ideals with a minimal
reduction $(y_1, \ldots, y_d)$ of $IJ$ obtained from a complete 
reduction of $I$ and $J$, then
$$
(\widetilde{I^aJ^b}) = \displaystyle{\bigcup_{k \geq 0}}I^{a+k}J^{b+k}
: (y_1^k, \ldots, y_d^k).
$$
\end{enumerate}
\el

\begin{proof}
(i). Let $x \in \widetilde{IJ}$, then $xI^nJ^n \subseteq I^{n+1}J^{n+1}$
for some $n$. Conversely if $xI^rJ^s
\subseteq I^{r+1}J^{s+1}$ for some $r,s \geq 0$ then for $n = max
\{r, s\}$, $xI^nJ^n \subseteq I^{n+1}J^{n+1}$ so that
$x \in (\widetilde{IJ})$.

(ii). By (i), $(\widetilde{I^aJ^b}) = \cup_{r,s \geq 0}I^{ar+a}J^{bs+b}
: I^{ar}J^{bs}$. 
Let $z \in (\widetilde{I^aJ^b})$
then for some $r, s$ we have $zI^{ar}J^{bs} \subseteq 
I^{ar+a}J^{bs+b}$. Set
$k = max\{ar, bs\}$. Then $zI^kJ^k \subseteq I^{a+k}J^{k+b}$ and hence
$z \in I^{a+k}J^{b+k} : I^kJ^k$. Let $zI^kJ^k \subseteq 
I^{a+k}J^{b+k}$ for some $k$. We may assume that $k
 = nab$ for $n \gg 0$. Therefore 
 $z \in I^{nab+a}J^{nab+b} : I^{nab}J^{nab} \subseteq (\widetilde{I^aJ^b})$. 

(iii). Suppose $z \in (\widetilde{I^aJ^b})$. Then for some $k$,
$zI^kJ^k \subseteq I^{a+k}J^{b+k}$, by (ii). Since $(y_1^k, \ldots,
y_d^k) \subseteq I^kJ^k$, we have $z(y_1^k, \ldots, y_d^k) \subseteq
I^{a+k}J^{b+k}$.  Let $zy_i^k \in I^{a+k}J^{b+k}$ for $i = 1,
\ldots, d$. Let $(\underline{y})$ denote the ideal $(y_1, \ldots,
y_d)$.  Then
$(IJ)^{m+n} = (\underline{y})^m(IJ)^n$ for all $m \geq 0$ and $n \geq
r_0 = r_{(\underline{y})}(IJ)$. Hence $(IJ)^{r+dk} =
(\underline{y})^{dk}I^rJ^r$ for $r \geq r_0.$  Therefore, 
$$
zI^{r+dk}J^{r+dk}  =  z(\underline{y})^{dk}I^rJ^r 
  =  \sum_{\sum i_j = dk}z y_1^{i_1}\cdots y_d^{i_d}I^rJ^r 
 \subseteq  I^{a+dk}J^{b+dk}I^rJ^r 
$$
Hence $z \in (\widetilde{I^aJ^b})$, by (ii). 
\end{proof}

\bl\label{artinrees}
 Let $I, J$ be ideals in a Noetherian ring $A$, $M$ a finite $A$-module
and $K$ an ideal of $A$ generated by $M$-regular elements. Then there exist
$t_1, t_2 > 0$ such that $I^rJ^sM :_M K =
I^{r-{t_1}}J^{s-{t_2}}(I^{t_1}J^{t_2}M :_M K)$ for all $r \geq t_1, s
\geq t_2$.
\el

\begin{proof}We follow the line of argument in [\cite{m}, Prop. 11.E].
Let $K = (a_1, a_2, \ldots, a_n)$ where $a_i$ are $M$-regular. Let $S$
be the multiplicatively closed subset generated by $a_1, \ldots, a_n$.
For $j = 1, \ldots, n$ consider the $A$-submodule $M_j = a_j^{-1}M$ 
of $S^{-1}M$ and set $L =
M_1 \oplus M_2 \oplus \cdots \oplus M_n$. Let $\Delta_M$ be the image
of the diagonal map $x \mapsto (\frac{x}{1}, \ldots, \frac{x}{1})$ from
$M$ to $L$. Since $a_i$'s are regular $\Delta_M \cong M$. Then, 
$$
I^rJ^sM :_M K = \bigcap_j \left(I^rJ^sM :_M a_j\right) = 
\bigcap_j (I^rJ^sM_j \cap M) \cong I^rJ^sL \cap \Delta_M.
$$ 
Since $L$ is a finite $A$-module and
$\Delta_M$ is a submodule of $L$, we can apply the generalized Artin-Rees
Lemma to get $t_1, t_2 > 0$ such that 
$$I^rJ^sL \cap \Delta_M =
I^{r-t_1}J^{s-t_2}(I^{t_1}J^{t_2}L \cap \Delta_M) \mbox{ for all } r \geq
t_1, s \geq t_2.
$$ 
Hence
$$ 
I^rJ^sM : K = I^{r-t_1}J^{s-t_2}(I^{t_1}J^{t_2}M : K) \mbox{ for
all } r \geq t_1, s \geq t_2.
$$
\end{proof}

\bl
Suppose $IJ$ has a reduction generated by regular elements, then for
$r, s \gg 0$, $(\widetilde{I^rJ^s}) = I^rJ^s$.
\el
\begin{proof} 

We first show that $I^{r+1}J^{s+1} : IJ = I^{r}J^{s}$ for $r,
s \gg 0$. Let $(\underline{\bf x}) = (x_1, \ldots, x_g)$ be a reduction
of $IJ$ generated by regular elements. Then, $I^nJ^n =
(\underline{\bf x})I^{n-1}J^{n-1}$ for $n \gg 0$ and hence $I^{r+1}J^{s+1} =
(\underline{\bf x})I^{r}J^{s}$ for $r, s \gg 0$. By setting $M = A$
and $K = (\underline{\bf x})$ in the Lemma \ref{artinrees}, we get $t_1, t_2 >
0$ such that $I^{r+1}J^{s+1} : (\underline{\bf x}) =
I^{r+1-t_1}J^{s+1-t_2}(I^{t_1}J^{t_2} : (\underline{\bf x})).$ Choose $r$
and $s$ large enough so that $r - t_1, s - t_2 \geq r_{(\underline{\bf
x})}(IJ).$ Then we have
\beqn
I^{r+1}J^{s+1} : IJ & \subseteq & I^{r+1}J^{s+1} : (\underline{\bf x}) \\
  & = & 
I^{r+1- t_1}J^{s+1-t_2}(I^{t_1}J^{t_2} : (\underline{\bf x})) \\
  & = &
(\underline{\bf x})I^{r- t_1}J^{s-t_2}(I^{t_1}J^{t_2} :
(\underline{\bf x})) \\
  & \subseteq &
I^{r}J^{s}.
\eeqn
Therefore $I^{r+1}J^{s+1} : IJ = I^{r}J^{s} ~\forall ~ r, s \gg 0$.
We claim that for all $k \geq 1$ and $r, s \gg 0$
$$ 
I^{r+k}J^{s+k} : I^kJ^k = I^rJ^s.
$$
Apply induction on $k$. The $k =1$ case has just been proved.
Let $k > 1$. Assume the result for $k-1$. Then
$$
I^{r+k}J^{s+k} : I^kJ^k  =  (I^{r+k}J^{s+k} : I^{k-1}J^{k-1}) : IJ = 
I^{r+1}J^{s+1} : IJ = 
I^rJ^s.
$$

\end{proof}

\section{A generalization of the Kirby-Mehran complex}

In this section we construct a bigraded analogue of a complex first
constructed by Kirby and Mehran in \cite{km}.
We study the cohomology modules of this complex and relate them to
those of the bigraded Rees algebras of two ideals.
Let $(A,\mm)$ be a d-dimensional Noetherian local ring with infinite
residue field and $I, J$ be
$\mm$-primary ideals of $A$. Let $\R$ and $\R^*$ be respectively 
the Rees and the extended Rees algebra of $A$ \wrt $I$ and $J$. 
Let $y_1, \ldots, y_m \in IJ$. For $k \geq 1$ set 
$(\underline{y})^{[k]} = (y_1^k, \ldots, y_n^k)$ and
$(\underline{yt})^{[k]} = ((y_1t_1t_2)^k, \ldots, (y_nt_1t_2)^k).$
Consider the Koszul complex $K^{\cdot}((\underline{yt})^{[k]}; \R)$ : 

$$
0 \lrar \R \lrar \R(k,k)^{n \choose 1} \lrar \cdots \lrar 
{\R((n-1)k, (n-1)k)}^{n \choose n-1}
\lrar \R(nk, nk) \lrar 0.
$$
This complex has a natural bigraded structure inherited from $\R$. 
Write the $(r,s)$th graded component, 
$K^{\cdot}_{(r,s)}((\underline{yt})^{[k]}; \R)$, of this complex :
$$
 0 \lrar (It_1)^r(Jt_2)^s \lrar {(It_1)^{r+k}(Jt_2)^{s+k}}^{n 
\choose 1} \lrar \cdots \lrar
(It_1)^{r+nk}(Jt_2)^{s+nk} \lrar 0.
$$

\noindent
This complex can be considered as a subcomplex of the Koszul complex :
$$
K^{\cdot}((\underline{y})^{[k]}; A):
0 \lrar A \lrar A^{n \choose 1} \lrar \cdots \lrar A^{n \choose n-1}
\lrar A \lrar 0
$$
Therefore there is map of complexes
$0 \lrar K^{\cdot}_{(r,s)}((\underline{yt})^{[k]}; \R) \lrar 
K^{\cdot}((\underline{y})^{[k]}; A).$ Since this inclusion is a chain map,
there exists a quotient complex. 

\begin{definition}
For $k \geq 1, r, s \in \ZZ$ and $n
\geq 1$ we define the complex $C^{\cdot}(n,k,r,s)$ to be the quotient
of the complex 
$K^{\cdot}((\underline{y})^{[k]}; A)$ by the complex
$K^{\cdot}_{(r,s)}((\underline{yt})^{[k]}; \R)$. 
\end{definition}

\noindent
We have the short exact sequence
\begin{eqnarray}\label{mehko}
0 & \lrar & K^{\cdot}_{(r,s)}((\underline{yt})^{[k]}; \R) \lrar 
K^{\cdot}((\underline{y})^{[k]}; A) \lrar C^{\cdot}(n,k,r,s) \lrar 0,
\end{eqnarray}

\noindent
One can easily see that  $C^{\cdot}(n,k,r,s)$ is the complex
$$
0 \lrar A/I^rJ^s \buildrel{d_C^0}\over\longrightarrow
\left(A/I^{r+k}J^{s+k}\right)^{{n}\choose{1}}
\buildrel{d_C^1}\over\longrightarrow \cdots
\buildrel{d_C^{n-1}}\over\longrightarrow
\left(A/I^{r+nk}J^{s+nk}\right) \buildrel{d_C^n}\over\lrar 0.
$$
where the differentials are induced by those of the Koszul complex
$K^{\cdot}(y_1^k, \ldots, y_n^k; A)$. We compute some of the
cohomology modules of this complex in the following proposition.

\bp\label{elimprop}
For all $k \geq 1, r,s \in \ZZ$  we have
\ben
\item[(i)] $H^0(C^{\cdot}(n,k,r,s)) = I^{r+k}J^{s+k} : (\underline{y}^{[k]})/I^rJ^s.$
\item[(ii)] $H^n(C^{\cdot}(n,k,r,s)) = A/(I^{r+k}J^{s+k} +
(\underline{y}^{[k]})).$ 
\item[(iii)]If $y_1, \ldots, y_n$ is an A-sequence, then
$$H^{n-1}(C^{\cdot}(n,k,r,s)) \cong \frac{(\underline{y}^{[k]}) \cap
I^{r+nk}J^{s+nk}}{(\underline{y}^{[k]})I^{r+(n-1)k}J^{s+(n-1)k}}.$$
\een
\ep

\bpf
\beqn
\mbox{(i)}~~ H^0(C^{\cdot}(n,k,r,s)) & = & \ker d^0_C \\
                            & = & \{\bar{u} \in A/I^rJ^s ~ \mid ~ y_i^ku \in
                                   I^{r+k}J^{s+k} \mbox{ for each } i = 1,
				   \ldots, n \} \\
                            & = & \frac{I^{r+k}J^{s+k} : (\underline{y})^{[k]}}{I^rJ^s}. \\
			    & ~ & \\
\mbox{(ii)}~~ H^n(C^{\cdot}(n,k,r,s)) & = & \frac{\ker d_C^n}{\im d_C^{n-1}} \\
                               & = &
			       \frac{A/I^{r+nk}J^{s+nk}}{(\underline{y})^{[k]} 
			       +I^{r+nk}J^{s+nk}/I^{r+nk}J^{s+nk}} \\
                               & \cong & \frac{A}{({\underline{y})^{[k]}
			       + I^{r+nk}J^{s+nk}}}.
\eeqn
(iii) Suppose that $y_1, \ldots, y_n$ is an A-sequence.
Consider the Koszul complex 

$$
K^{\cdot}((\underline{y})^{[k]}, A) :  \hspace*{.3in}\cdots \lrar
A^{n \choose {n-2}} \buildrel{d^{n-2}_K}\over\longrightarrow A^{n
\choose {n-1}} \buildrel{d^{n-1}_K}\over\longrightarrow (y_1^k,
\ldots, y_n^k) \lrar 0
$$
Since $(y_1^k, \ldots, y_n^k)$ is an $A$-sequence, this is an exact
sequence. Tensoring by \\ $A/I^{r+{(n-1)}k}J^{s+{(n-1)}k}$, we get an exact
sequence
$$
\left(\frac{A}{I^{r+{(n-1)}k}J^{s+{(n-1)}k}}\right)^{n \choose {n-2}}
\buildrel{\bar{d}^{n-2}_K}\over\longrightarrow
\left(\frac{A}{I^{r+{(n-1)}k}J^{s+{(n-1)}k}}\right)^{n \choose {n-1}}
$$
~
$$
\buildrel{\bar{d}^{n-1}_K}\over\longrightarrow
\frac{(\underline{y})^{[k]}}{(\underline{y})^{[k]}
I^{r+{(n-1)}k}J^{s+{(n-1)}k}} \lrar 0.
$$

We have $\im \bar{d}^{n-2}_K = \im d^{n-2}_C$ and a
commutative diagram of exact rows

$$
\begin{CD}
0 @>>> \im \bar{d}^{n-2}_K @>>> 
\left(\frac{A}{I^{r+{(n-1)}k}J^{s+{(n-1)}k}}\right)^n @>>>
\frac{(y_1^k, \ldots, y_n^k)}{(y_1^k, \ldots,
y_n^k)I^{r+{(n-1)}k}J^{s+{(n-1)}k}} @>>> 0\\
& & @V{\alpha}VV @V{id}VV @V{\gamma}VV \\
0 @>>> \ker d^{n-1}_C @>>>
\left(\frac{A}{I^{r+{(n-1)}k}J^{s+{(n-1)}k}}\right)^n @>>>
\frac{A}{I^{r+nk}J^{s+nk}} 
\end{CD}
$$
where $\alpha$ is the inclusion map and $\gamma$ is the natural map.
By the Snake lemma, we get 
$$
H^{n-1}(C^{\cdot}(n,k,r,s)) \cong coker \alpha \cong \ker\gamma\cong
\frac{(y_1^k, \ldots, y_n^k) \cap I^{r+nk}J^{s+nk}}{(y_1^k, \ldots,
y_n^k)I^{r+(n-1)k}J^{s+(n-1)k}}.
$$
\epf

For the rest of the section let $I$ and $J$ be $\mm$-primary ideals
of $A$. Let $x_{1j} \in I$ and $x_{2j} \in J$ for $j = 1, \ldots, d$
and set $y_i = x_{1i}x_{2i}.$ 

\bp\label{reescoho}
Let $r,s \in \ZZ$.
\ben
\item[(i)] For all $k \geq 1$, there is an exact sequence of $A$-modules
$$
0  \lrar   H^0((\underline{yt})^{[k]}; \R)_{(r,s)} 
\lrar H^0((\underline{y})^{[k]}; A) \lrar H^0(C^{\cdot}(n,k,r,s)) 
\lrar H^1((\underline{yt})^{[k]}; \R)_{(r,s)} \lrar \cdots
$$
\item[(ii)] There is an exact sequence of $A$-modules
$$
0 \lrar H^0_{(\underline{yt})}(\R)_{(r,s)} \lrar
H^0_{(\underline{y})}(A) \lrar
\displaystyle\lim_{\stackrel{\longrightarrow}{k}} 
H^0(C^{\cdot}(n,k,r,s))
\lrar H^1_{(\underline{yt})}(\R)_{(r,s)} \lrar \cdots
$$
\een
\ep

\bpf
(i). Follows from the long exact sequence of Koszul homology modules 
corresponding to (\ref{mehko}).

(ii). For each $i$, consider the commutative diagram of complexes
$$
\begin{CD}
K^{\cdot}((y_it_1t_2)^k; \R):  0  @>>> \R @>{(y_it_1t_2)^k}>> \R @>>> 0 \\ 
& &  @V{id}VV  @V{y_it_1t_2}VV   \\
K^{\cdot}((y_it_1t_2)^{k+1}; \R):  0  @>>>  \R  @>{(y_it_1t_2)^{k+1}}>> 
\R  @>>>  0. 
\end{CD}
$$

This gives a map $\otimes_{i = 1}^n K^{\cdot}((y_it_1t_2)^k; \R) \lrar 
\otimes_{i = 1}^n K^{\cdot}((y_it_1t_2)^{k+1}; \R)$, i.e., we get a map 

$$
K^{\cdot}((\underline{yt})^{[k]}; \R) \lrar
K^{\cdot}((\underline{yt})^{[k+1]}; \R)
$$
and its restriction to the $(r,s)$-th component gives the map 

$$
K^{\cdot}_{(r,s)}((\underline{yt})^{[k]}; \R) \lrar
K^{\cdot}_{(r,s)}((\underline{yt})^{[k+1]}; \R). 
$$
Thus we obtain a commutative diagram of exact sequences

$$
\begin{CD}
0  @>>>  K_{(r,s)}((\underline{yt})^{[k]}; \R) @>>>
K^{\cdot}((\underline{y})^{[k]}; A) @>>> C^{\cdot}(n,k,r,s) @>>> 0\\
& &  @VVV  @VVV @VVV  \\
0  @>>>  K_{(r,s)}((\underline{yt})^{[k+1]}; \R) @>>>
K^{\cdot}((\underline{y})^{[k+1]}; A) @>>> C^{\cdot}(n,k+1,r,s) @>>> 0\\
\end{CD}
$$
Apply $\displaystyle\lim_{\stackrel{\longrightarrow}{k}}$
to the long exact sequence of the cohomology modules to get (ii).
\epf

\bc\label{cor1}
Let $(A,\mm)$ be Cohen-Macaulay of dimension $d \geq 2$ and  
$(x_{ij})$; $i = 1,2$; $1 \leq j \leq d$ be a 
complete reduction of $(I, J).$ Let $r, s \in \ZZ$. Then
\ben
\item[(i)] For all $k \geq 0$, we have 
$$
H^i((\underline{yt})^{[k]}; \R)_{(r,s)} \cong
H^{i-1}(C^{\cdot}(d,k,r,s)) ~~~ \fa~~ 0 \leq i \leq d-1.
$$
and an exact sequence of $A$-modules
$$
0  \lrar  H^{d-1}(C^{\cdot}(d,k,r,s)) \lrar 
H^{d}((\underline{yt})^{[k]}; (\R))_{(r,s)}  \lrar
H^d((\underline{y})^{[k]}; A) 
\lrar H^{d}(C^{\cdot}(d,k,r,s)) \lrar 0
$$
\item[(ii)] There is an isomorphism of $A$-modules

$$
H^i_{(\underline{yt})}(\R)_{(r,s)} \cong 
\displaystyle\lim_{\stackrel{\longrightarrow}{k}}
H^{i-1}(C^{\cdot}(d,k,r,s)) ~~~ \fa~~ 0 \leq i \leq d-1
$$
and an exact sequence 
\beqn
0 & \lrar & \displaystyle\lim_{\stackrel{\longrightarrow}{k}}
\frac{(\underline{y})^{[k]} \cap I^{r+dk}J^{s+dk}}
{(\underline{y})^{[k]} I^{r+(d-1)k}J^{s+(d-1)k}} \lrar
H^d_{(\underline{yt})}(\R)_{(r,s)} \lrar 
H^d_m(A) \\
 & \lrar & \displaystyle\lim_{\stackrel{\longrightarrow}{k}}
\frac{A}{(\underline{y})^{[k]} + I^{r+dk}J^{s+dk}} \lrar 0.
\eeqn
\item[(iii)] 
$
H^1_{(\underline{yt})}(\R)_{(r,s)}~~ \cong~~ 
\frac{(\widetilde{I^rJ^s})}{I^rJ^s}. 
$
\een
\ec
\bpf
(i) Consider the long exact sequence of cohomology modules
corresponding to (\ref{mehko}). 
\begin{eqnarray*}
0 & \lrar & H^0(K^{.}((\underline{yt})^{[k]}; \R)) \lrar
H^0(K^.((\underline{y})^{[k]}; A)) \lrar H^0(C^.(d,k,r,s)) \\
& \lrar & H^1(K^{.}((\underline{yt})^{[k]}; \R)) \lrar \cdots
\end{eqnarray*}
Since $A$ is Cohen-Macaulay $H^i(K^.((\underline{y})^{[k]}; A) = 0$
for all $0 \leq i \leq d-1$. Hence (i) follows.

\noindent
(ii) Apply $\displaystyle\lim_{\stackrel{\longrightarrow}{k}}$  
to (i).

\noindent
(iii) By (ii) and Lemma \ref{rrprod} we have 
$$
H^1_{(\underline{yt})}(\R)_{(r,s)}  \cong  
\displaystyle\lim_{\stackrel{\longrightarrow}{k}}H^0(C^{\cdot}(d,k,r,s))
 =  \displaystyle\lim_{\stackrel{\longrightarrow}{k}}
\frac{I^{r+k}J^{s+k} : (\underline{y})^{[k]}}{I^rJ^s} 
 =  \frac{(\widetilde{I^rJ^s})}{I^rJ^s}.
$$
\epf
A similar theory can be developed for the extended Rees algebra 
by setting $I^r = A = J^s$ if $r, s \leq 0$ and
defining the complex $C^{\cdot}(n,k,r,s)^*$ in a similar way as we
defined $C^{\cdot}(n,k,r,s)$. We can prove results similar to
Proposition \ref{elimprop}, Proposition \ref{reescoho} etc. 
First we prove a general result relating local cohomology modules of
two bigraded algebras which will help us in relating 
the local cohomology modules of the Rees and the extended Rees algebras.

\bp\label{locrr*}
Let $R = \oplus_{r,s \geq 0}R_{(r,s)} \hookrightarrow \oplus_{r,s \in
\ZZ}R_{(r,s)} = R^*$ be an inclusion of bigraded algebras over $R_{(0,0)}$,
a Noetherian ring. Then 
\ben
\item[(i)] For $i > 1$, we have $H^i_{R_{++}}(R) \cong H^i_{R_{++}}(R^*)$.
\item[(ii)] We have an exact sequence 
$$
0 \lrar H^0_{R_{++}}(R) \lrar H^0_{R_{++}}(R^*) \lrar R^*/R \lrar
H^1_{R_{++}}(R) \lrar H^1_{R_{++}}(R^*) \lrar 0.
$$
\een
\ep

\begin{proof}
Consider the exact sequence of bigraded $R$-modules.
\begin{eqnarray}\label{sesrr*1}
0 \lrar R \lrar R^* \lrar R^*/R \lrar 0.
\end{eqnarray}

Since $R_{++}$ acts nilpotently on $R^*/R$, 
$ H^0_{R_{++}}(R^*/R)   =  R^*/R$ and $H^i_{R_{++}}(R^*/R) = 0$ for all
$i \neq 0.$
The proposition follows from the long exact sequence of local
cohomology modules derived from (\ref{sesrr*1}).

\end{proof}

\bc\label{gen1}
Consider the bigraded rings $\R = A[It_1, Jt_2] \hookrightarrow 
\R^*= A[It_1, Jt_2, t_1^{-1}, t_2^{-1}]$ and
$\G = \oplus_{r,s \geq 0}I^rJ^s/I^{r+1}J^{s+1} \hookrightarrow
\G^* = \R^*/t_1^{-1}t_2^{-1}\R^*$. Then
\ben
\item[(i)]For all $i \geq 2$ we have the isomorphism
$ H^i_{\R_{++}}(\R) \cong H^i_{\R_{++}}(\R^*)$ and 
there is an exact sequence of bigraded $\R$-modules
$$
0 \lrar H^0_{\R_{++}}(\R) \lrar  H^0_{\R_{++}}(\R^*)
\lrar \R^*/\R \lrar H^1_{\R_{++}}(\R) \lrar
H^1_{\R_{++}}(\R^*) \lrar 0.
$$

\item[(ii)]For all $i \geq 2$ we have $H^i_{\G_{++}}(\G) \cong
H^i_{\G_{++}}(\G^*)$ and
there is an exact sequence of bigraded $\G$-modules
$$
0 \lrar H^0_{\G_{++}}(\G) \lrar  H^0_{\G_{++}}(\G^*)
\lrar \G^*/\G \lrar H^1_{\G_{++}}(\G) \lrar
H^1_{\G_{++}}(\G^*) \lrar 0.
$$
\een
\ec

\bc\label{cor2} For all $r, s \geq 0,$
$$
H^1_{\R_{++}}(\R^*)_{(r,s)}  \cong  
\frac{(\widetilde{I^rJ^s})}{I^rJ^s}. 
$$
\ec
\begin{proof}
Use Corollary \ref{cor1}(iii) and Corollary \ref{gen1}(i) to get the
required result.
\end{proof}

\section{The difference formula}

In this section we obtain an expression for the difference of
Bhattacharya polynomial and Bhattacharya function. The main motivation
were  results of Johnston-Verma \cite{jv} and C. Blancafort \cite{bl} 
which express the difference of
Hilbert-Samuel polynomial and Hilbert-Samuel function in terms of
the Euler characteristic of the   Rees algebra (resp. extended Rees algebra).
We have followed Blancafort's elegant line
of approach in the proof. However, we prove the theorem only for
non-negative integers. The question remains still open for negative
integers. 

\begin{theorem}\label{main} 
Let $\R^* = A[It_1, Jt_2, t_1^{-1}, t_2^{-1}]$. Then 

\begin{enumerate}
\item[(i)] $\lm_A(H^i_{\R_{++}}(\R^*)_{(r,s)}) < \infty$ 
$\fa~ r, s \in \ZZ; i = 0, 1, \ldots, d.$
\item[(ii)] $P(r, s) - B(r, s) = \sum_{i = 0}^d (-1)^i
\lm_A(H^i_{\R_{++}}(\R^*)_{(r,s)}) \fa ~ r, s \geq 0.$
\end{enumerate}
\end{theorem}

\begin{proof}(i). By  Theorem \ref{imp2}, $H^i_{\R_{++}}(\R)_{(r,s)}$ 
are finitely generated $A$-modules and they vanish for $r,s \gg 0$. 
By Lemma \ref{imp1} and Corollary \ref{gen1},
$H^i_{\R_{++}}(\R^*)_{(r,s)} = 0$ for all $r,s \gg 0$. 
We have an exact sequence of bigraded $\R$-modules :
\begin{eqnarray}\label{rrg}
0 \lrar \R^*(1,1) \buildrel{t_1^{-1}t_2^{-1}}\over\longrightarrow \R^*
\lrar \G^* \lrar 0,
\end{eqnarray}
where $\G^* = \R^*/t_1^{-1}t_2^{-1}\R^*.$ 
By the change of ring principle,
$H^i_{\R_{++}}(\G^*) = H^i_{\G_{++}}(\G^*)$ for
all $i \geq 0$. 
From the above short exact sequence  we obtain the long exact sequence :
$$
0 \lrar H^0_{\R_{++}}(\R^*)_{(r+1, s+1)} \lrar
H^0_{\R_{++}}(\R^*)_{(r, s)} \lrar H^0_{\G_{++}}(\G^*)_{(r, s)} 
\lrar H^1_{\R_{++}}(\R^*)_{(r+1, s+1)} \lrar \cdots
$$
We prove (i) by decreasing induction on $r$ and $s$.
Since $H^i_{\R_{++}}(\R^*)_{(r, s)} = 0$ for all $r,s \gg 0$, the
result is obviously true for $r,s \gg 0$. 
Consider the exact sequence
$$
\cdots \lrar H^i_{\R_{++}}(\R^*)_{(r+1, s+1)} \lrar 
H^i_{\R_{++}}(\R^*)_{(r, s)} \lrar H^i_{\G_{++}}(\G^*)_{(r, s)} \lrar
\cdots
$$
By induction $H^i_{\R_{++}}(\R^*)_{(r+1, s+1)}$ has finite length. By
Theorem \ref{imp2} and Corollary \ref{gen1}(ii) 
$H^i_{\G_{++}}(\G^*)_{(r, s)}$ is a finitely generated $\G_{00}$-module. 
Since $\G_{00}$ is Artinian $H^i_{\G_{++}}(\G^*)_{(r, s)}$ has finite 
length. Therefore $H^i_{\R_{++}}(\R^*)_{(r, s)}$ has finite length. \\
(ii). For a bigraded module $M$ over the bigraded ring $\R$, set
$$\chi_M(r,s)  =  \sum_{i \geq 0}(-1)^i
\lambda_A(H^i_{\R_{++}}(M)_{(r, s)}) ~~ \and~~ g(r,s)  =  P(r,s) - B(r,s).$$
Then from the exact sequence (\ref{rrg}) we get for all $r,s \geq 0,$
\begin{eqnarray*}
\chi_{\R^*(1,1)}(r,s) - \chi_{\R^*}(r,s) & = & \chi_{\R^*}(r+1, s+1) - 
\chi_{\R^*}(r,s) \\
& = & - \chi_{\G^*}(r,s)  = - \chi_{\G}(r,s)\hspace*{0.3in} (\mbox{by  4.6(ii)})\\
& = & P_{\G}(r,s) - H_{\G}(r,s) = P_{\G^*}(r,s) - H_{\G^*}(r,s) \\
& = & \left( P(r+1,s+1) - P(r,s) \right) - \left(B(r+1, s+1) - B(r,s) 
      \right) \\
& = & g(r+1, s+1) - g(r,s).
\end{eqnarray*}
Set $h(r,s) = \chi_{\R^*}(r,s)  - g(r,s)$. Then $h(r,s) = h(r-1, s-1)$ 
for all $r,s \geq 0$ and $h(r,s) = 0$ for all $r,s \gg 0.$ This 
clearly implies that $h(r,s) = 0$ for all $r,s \geq 0.$
\end{proof}

\begin{corollary}\label{app1}
Let $(A,\mm)$ be a $2$-dimensional Cohen-Macaulay local ring and $I, J$ be
$\mm$-primary 
ideals of $A$. Then for all $r, s \geq 0$
$$ P(r, s) - B(r, s) = \lm(H^2_{\R_{++}}(\R)_{(r,s)}) -
\lm(\widetilde{I^rJ^s}/I^rJ^s).$$
In particular $$ e_{00} = \lm(H^2_{\R_{++}}(\R)_{(0,0)}).$$
\end{corollary}

\bpf
By the previous theorem, 
$$
P(r,s) - B(r,s) = \lm(H^0_{\R_{++}}(\R)_{(r,s)}) - 
\lm(H^1_{\R_{++}}(\R)_{(r,s)}) + \lm(H^2_{\R_{++}}(\R)_{(r,s)}).
$$
Since $I$ and $J$ are $\mm$-primary, $\R_{++}$ contains a regular
element. Therefore $H^0_{\R_{++}}(\R) = 0.$ By Proposition \ref{gen1}, 
$$
H^1_{\R_{++}}(\R)_{(r,s)} \cong \frac{\widetilde{I^rJ^s}}{I^rJ^s}.
$$
Now, 
$$ e_{00} = P(0, 0) - B(0,0)= \lm(H^2_{\R_{++}}(\R)_{(0,0)}).$$

\epf

\section{Bigraded Cohen-Macaulay Rees Algebras }

In the previous section we have established  a formula for the difference
between
the Bhattacharya function and Bhattacharya polynomial. It is
interesting to know when is the
Bhattacharya function equal to the Bhattacharya polynomial. 
Here we
give a partial answer to this question, in dimension 2. Huneke
(Theorem 2.1, \cite{h2}) and
Ooishi (Theorem 3.3, \cite{o}) gave a characterization for the 
reduction number of an $\mm$-primary ideal to be at most 1 in terms of
$e_0(I)$ and $e_1(I)$. Huckaba and Marley (Corollary 4.8, Corollary
4.10, \cite{hm})
generalized this result for higher reduction numbers. In particular,
they characterized Cohen-Macaulay property of the Rees algebra to be
Cohen-Macaulay in terms of the $e_1(I)$. It is natural to ask whether
one can characterize the Cohen-Macaulay property of bigraded Rees
algebras in terms of coefficients of the Bhattacharya polynomial. The
Theorem \ref{jr} below answers this in dimension 2. A similar
characterization for Cohen-Macaulayness of the multi-Rees algebras in
higher dimension in terms of Bhattacharya coefficients is not known.

We need another generalization of reductions for two ideals, namely 
joint reductions. 
Let $A$ be a commutative ring with identity and let $I_1, I_2, \ldots,
I_g$ be ideals of $A$. A system of elements $(\underline{\bf x}) :=
(x_1, x_2, \ldots, x_g)$, where $x_i \in I_i$, is said to be a {\it joint
reduction} of the sequence of ideals $(I_1, I_2, \ldots, I_g)$
if there exist positive integers  $d_1, d_2, \ldots, d_g$ such
that 
$$ x_1I_1^{d_1-1}I_2^{d_2}\cdots I_g^{d_g} + \cdots +
x_gI_1^{d_1}\cdots I_{g-1}^{d_{g-1}}I_g^{d_{g}-1} = I_1^{d_1}\cdots
I_g^{d_g}.$$
We say that the sequence of ideals $(I_1, \ldots, I_g)$ has {\it joint
reduction number zero } if
$$x_1I_2\cdots I_g + \cdots + x_gI_1\cdots I_{g-1} = I_1I_2\cdots I_g.$$
We first prove a general property of the Bhattacharya coefficients.

\begin{lemma}\label{e10d1}
Let $(A,\mm)$ be a $1$-dimensional Cohen-Macaulay local ring with
infinite residue field. Let $I$ and $J$ be $\mm$-primary ideals of
$A$. Then
\begin{enumerate}
\item[(i)] $P(r+1, s) - H(r+1,s) \geq P(r,s) - H(r,s)$  and
\item[] $P(r, s+1) - H(r,s+1) \geq P(r,s) - H(r,s)$.
\item[(ii)] $\lambda(A/I) \geq e_{10} + e_{00}$ and
$\lambda(A/J) \geq e_{01} + e_{00}$.
\end{enumerate}
\end{lemma}

\begin{proof}
Let
$(x) \subseteq I$ be a reduction of $I$. Then
\begin{eqnarray*}
P(r+1, s) - H(r+1, s) & = & e_{10}(r+1)+e_{01}s+e_{00} - 
\lambda(A/I^{r+1}J^s)\\
 & = & P(r,s) + e_{10} - \lambda(A/I^{r+1}J^s) \\
 & \geq & P(r,s) + \lambda(A/(x)) - \lambda(A/xI^rJ^s) \\
 & = & P(r,s) - \lambda((x)/xI^rJ^s) \\
 & = & P(r,s) - H(r,s)
\end{eqnarray*}

Similarly one can prove that $P(r, s+1) - H(r,s+1) \geq P(r,s) - 
H(r,s)$. 
From (i) it is clear that $P(r,s) - H(r,s) \leq 0$ for all $r,s$.
Putting $(r,s) = (1,0)$ and $(r,s) = (0, 1)$ we get  (ii).
\end{proof}

\begin{lemma}\label{e10d2}
Let $(A,\mm)$ be a $2$-dimensional Cohen-Macaulay local ring and $I$,
$J$ be $\mm$-primary ideals of $A$. Then $\lambda(A/I) \geq e_{10}$ 
\mbox{and} $\lambda(A/J) \geq e_{01}$.
\end{lemma}

\begin{proof}
Let $(x,y)$, where $x \in I$ and $y \in J$,  be a joint reduction of 
$(I,J)$. Choose the joint
reduction such that $x$ is superficial for $I$ and $J$. Let
$~\bar{}~$ denote ``modulo x". Let $\bar{H}(r,s)$ and $\bar{P}(r,s)$ denote
the Bhattacharya function and Bhattacharya polynomial of the 
$\bar{\mm}$-primary ideals $\bar{I}$ and $\bar{J}$ of $\bar{A} =
A/(x)$. \\
{\it Claim :} $\bar{P}(r,s) = P(r,s) - P(r-1, s)$. \\
From the following exact sequence
$$
\begin{CD}
0 @>>> I^rJ^s : x/I^rJ^s @>>> A/I^rJ^s @>x >> A/I^rJ^s
@>>> A/(I^rJ^s, x) @>>> 0
\end{CD}
$$
$\lambda(I^rJ^s : x/I^rJ^s) = \lambda (A/I^rJ^s, x)$. Then for all
$r,s \gg 0$,
\begin{eqnarray*}
\bar{P}(r,s) & = & \lambda(A/\bar{I}^r\bar{J}^s) =\lambda(A/(I^rJ^s,
x)) \\
 & = & \lambda(I^rJ^s : x/I^rJ^s) \\
 & = & \lambda(I^{r-1}J^s/I^rJ^s)
 \hspace*{0.2in}(\mbox{since x is superficial for I and J}) \\
 & = & P(r,s) - P(r-1,s)
\end{eqnarray*}
Therefore 
\begin{eqnarray*}
\bar{P}(r,s) & = & e_{20} \left[{r \choose 2} - {r-1 \choose 2}
\right] + e_{11}(r - (r-1))s + e_{10}(r - (r-1)) \\
& = & e_{20}(r-1) + e_{11}s + e_{10} \\
& = & e_{20}r + e_{11}s + e_{10} - e_{20}.
\end{eqnarray*}
Since $\dim \bar{A} = 1$, by Lemma \ref{e10d1},
$\lambda(\bar{A}/\bar{I}) \geq e_{20} + (e_{10} - e_{20})$. Hence
$\lambda(A/I) \geq e_{10}$. Similarly one can prove that $\lambda(A/J)
\geq e_{01}$.
\end{proof}

\begin{theorem}\label{jr}
Let $(A,\mm)$ be a $2$-dimensional Cohen-Macaulay local ring and $I, J$ be
$\mm$-primary ideals of $A$. Let $P(r, s) = \sum_{i+j \leq 2}
e_{ij}{r \choose i}{s \choose j}$ be the Bhattacharya polynomial of
$I$ and $J$ corresponding to the function $B(r, s) =
\lm(A/I^rJ^s)$. Then the following conditions are equivalent:
\begin{enumerate}
\item[(1)] $e_{10} = \lm(A/I)$ and $e_{01} = \lm(A/J).$
\item[(1$^\prime$)] $e_{10} \geq \lm(A/I)$ and $e_{01} \geq \lm(A/J)$.
\item[(2)] $P(r, s) = B(r,s)$ for all $r, s \geq 0.$ 
\item[(3)] The joint reduction number of $(I, J)$ is zero, $r(I) \leq
1$ and $r(J) \leq 1$.
\item[(4)] The Rees ring $A[It_1, Jt_2]$ is Cohen-Macaulay.
\end{enumerate}
\end{theorem}

\begin{proof}

The equivalence of (1) and (1$^\prime$) is clear from Lemma \ref{e10d2}.
First we show that hypotheses in (1) imply that the joint reduction 
number of $(I,J)$ is zero.  
By Theorem 3.2 of \cite{v}, it is enough to show that $e_1(I|J) = 
\lm(A/IJ) - \lm(A/I) -
\lm(A/J)$. By
Corollary \ref{app1}
\begin{eqnarray*}
e_{00} & = & \lm(H^2_{\R_{++}}(\R)_{(0,0)}) \\
e_1(I|J) + e_{10} + e_{01} + e_{00} - \lm(A/IJ) &= &
\lm(H^2_{\R_{++}}(\R)_{(1,1)}) - \lm(\widetilde{IJ}/IJ).
\end{eqnarray*}
Let $(y_1, y_2)$ be a reduction of $IJ$ coming from a complete
reduction of $(I, J)$. It follows from the long exact sequence 
of local cohomology modules
corresponding to the short exact sequence 
$$
0 \lrar \R^*(-1, -1)  \buildrel{.y_1t_1t_2}\over\longrightarrow\R^* \lrar
\R^*/y_1t_1t_2\R^* \lrar 0
$$
and Proposition \ref{locrr*}, that for all $r, s \in \ZZ$
$$\lm(H^2_{\R_{++}}(\R)_{(r+1,s+1)}) \leq
\lm(H^2_{\R_{++}}(\R)_{(r,s)}).$$
Therefore 
 $$e_1(I|J) + e_{10} + e_{01} + e_{00} - \lm(A/\widetilde{IJ}) \leq
 e_{00}.$$ 
Hence 
\begin{eqnarray*}
e_1(I|J) & \leq & \lm(A/\widetilde{IJ}) - \lm(A/I) - \lm(A/J) \\
& \leq & \lm(A/IJ) - \lm(A/I) - \lm(A/J).
\end{eqnarray*}
By the isomorphism $A/I \oplus A/J \cong (a,b)/aJ+bI$ for any
regular sequence
$(a, b)$ where $a \in I, \mbox{ and } b \in J$, it follows that
$$e_1(I|J) \geq \lm(A/IJ) - \lm(A/I) - \lm(A/J).$$
Therefore 
$$e_1(I|J) = \lm(A/IJ) - \lm(A/I) - \lm(A/J).$$
Since the joint reduction number of $(I, J)$ is zero, by Theorem 3.2 of 
\cite{v}, for all $r, s \geq 1$
$$\lm(A/I^rJ^s) = \lm(A/I^r) + e_1(I|J) rs + \lm(A/J^s).$$
$$\mbox{Write }\lm(A/I^r) = e_0(I){r \choose 2} + e_1(I) r + e_2(I) \mbox{ and }
\lm(A/J^s) =  e_0(J){s \choose 2} + e_1(J) s + e_2(J).$$
The reader may note that this way of writing the Hilbert polynomials 
of $I$ and $J$ is different from the way in which the Hilbert
polynomial is usually written. Therefore the first Hilbert
coefficient $e_1(I)$ appearing in the formulas above is different from
the $e_1(I)$ appearing in papers of, for example, Huneke and Ooishi.
Therefore, for $r, s \gg 0,$ we have, 
$$
P(r,s ) =  e_0(I){r \choose 2} + e_1(I|J) rs + e_0(J){s \choose 2} 
+ e_1(I) r +  e_1(J) s + e_2(I) + e_2(J).
$$
By assumption $e_1(I) = \lm(A/I)$ and $e_1(J) = \lm(A/J)$. By the
Huneke-Ooishi theorem, \cite{h2}, for $d =2$ we have
$r(I) \leq 1$, $e_2(I) = 0$ and
$r(J) \leq 1$, $e_2(J) = 0$. This proves (3) as well as (2).
The statement (2) $\Rightarrow$ (1) is obvious. The equivalence of (2)
and (3) follows from Theorem 3.2 of \cite{v} and Theorem 2.1 of
\cite{h2}.
The equivalence of (3) and (4) follows from Corollary 3.5 of \cite{hy}
and Goto-Shimoda Theorem \cite{gs}.
\end{proof}
The following example shows that a naive generalization of Theorem
\ref{jr} does not work for $d > 2$.

\begin{example}
Let $A = k[\![x,y,z]\!]$, $I = (x^2, xy, y^2, z)$ and $J = (x,y^3,z)$.
Then $(x^2, y^2, z)$ is a reduction of $I$ with reduction number 1.
One can also check that $IJ = (x,z)I + y^2J = xI + (y^2,z)J$.
Therefore $r(I) = 1$, $r(J) = 0$ and joint reduction numbers of $(I,
J)$ are zero.  One can see from computations on Macaulay 2 \cite{gs2} 
that $\depth \R = 4$.  But $\dim \R = 5$. Therefore $\R$
is not Cohen-Macaulay.
\end{example}

\begin{example}
Consider the plane curve $f = y^2 - x^n = 0.$ Put $A = \CC[[x, y]]$
and $\mm = (x,y)A.$ Let $J$ denote the Jacobian ideal $(f_x, f_y)$ of
$f = 0.$ Then $r(J) = r(\mm) = 0.$ Moreover, $y\mm+xJ = \mm J.$
Therefore by the previous theorem, the Bhattacharya polynomial of
$\mm$ and $J$ is given by the formula 
$$\lm(A/\mm^rJ^s) = {r \choose 2} + rs + (n-1){s \choose
2} + r + (n-1) s \mbox{ for all } r, s \geq 0.$$
\end{example}

\begin{example}
We give an example to show that neither of the conditions in (1) of
Theorem \ref{jr} can be dropped to get the conclusions (2) and (3).
Let $(A,\mm)$ denote a 2-dimensional regular local ring. Let $\mm = (x,
y)$ and $I = (x^3, x^2y^4, xy^5, y^7)$. Then $I\mm = x^3\mm + yI$. By
Theorem 3.2 of \cite{v}, we get 
\begin{eqnarray*}
\lambda(A/\mm^rI^s) & = & \lambda(A/\mm^r) + e_1(\mm|I) rs + 
                          \lambda(A/I^s)\\
                  & = & {r+1 \choose 2} + o(I) rs + \lambda(A/I^s).
\end{eqnarray*}
In the above equation $o(I)$ denotes the $\mm$-adic order of $I$ which
is 3. The fact that $e_1(\mm|I) = o(I)$ is proved in \cite{v}. We now
calculate the Hilbert polynomial of $I$.

The ideal $J = (x^3, y^7)$ is a minimal reduction of  $I$ and 
$JI^2 = I^3$ and $\lambda(I^2/JI) = 1.$ By a result of Sally,
\cite{s3}, $\lambda(R/I^n) = P_I(n)$ for all $n > 1$. Here $P_I(n)$
denotes the Hilbert polynomial of $I$ corresponding to the Hilbert
function $\lambda(A/I^n)$. By using Macaulay 2 \cite{gs2}, we find that
$\lambda(A/I) = 16, \lambda(A/I^2) = 52, \lambda(A/I^3) = 109$. 
Therefore the Hilbert
polynomial $P_I(n) = 21{n+1 \choose 2} - 6{n \choose 1} + 1$.
Hence the Bhattacharya polynomial is
\beqn
P(r,s)& = &{r+1 \choose 2} + 3 rs + 21{s+1 \choose 2} - 6{s \choose 1} 
	     + 1 \\
      & = &{r \choose 2} + 3rs + 21{s \choose 2} + {r \choose 1} +
            15{s \choose 1} + 1.
\eeqn
Therefore $e_{01} = 15 < \lambda(R/I).$ Notice that the constant term
of the Bhattacharya polynomial is non-zero. 
\end{example}

\end{document}